\documentclass [a4paper,10pt]{article}
\usepackage [utf8] {inputenc}
\usepackage {amsmath}
\usepackage {amsthm}
\usepackage {amssymb}
\usepackage {amsfonts}
\newtheorem {lemma}{Lemma}
\newtheorem {theorem}{Theorem}

\title {The distribution of $d_4(n)$ in arithmetic progressions}
\author {Tomos Parry}
\date {}
\begin {document}
\maketitle 
\begin {abstract}
We use the Petrow-Young \cite {petrowyoung} subconvexity bound for Dirichlet $L$-functions to show that $d_4(n)$ has exponent of distribution $4/7$ when we allow an average over $a$ mod $q$, thereby giving an equidistribution result for $d_4(n)$ which goes past the $1/2$ barrier for the first time.
\end {abstract}
\begin {center}
\section {Introduction}
\end {center}
Let $d_k(n)$ be the $k$-fold divisor function
\[ d_k(n)=\sum _{d_1\cdot \cdot \cdot d_k=n}1\]
and write $d(n)=d_2(n)$ for the classical divisor function. What we would like is equidistribution over the residue classes in the form 
\[ \sum _{n\leq x\atop {n\equiv a(q)}}d_k(n)=\text { main term }+\text { small error term }\]
where the main term is a simple function of size around $x/q$ and where the error term is a power-saving on the main term - important is the range of validity of $q$ and in particular whether we can take $q$ past $\sqrt x$.
\\
\\ For $k=2$ we have a classical result of Selberg and of Hooley that says we have equidistribution for $q$ up to around $x^{2/3}$ and this has resisted much improvement. For $k=3$ the  work of Friedlander and Iwaniec \cite {fi} in the 80's and then improvements by Fouvry, Kowalski and Michel \cite {kowalski} in 2015 meant we can take prime $q$ up to $x^{1/2+1/46}$ for equidistribution. For $k\geq 4$ we don't know if we can take $q$ past $\sqrt x$, although for $k=4$ we can take it to be exactly $\sqrt x$ from the work of Linnik. 
\\
\\ Our interest in equidistribution results for higher divisor functions comes from the links to bounded gaps between primes - see, for example, Remark 1.4 of \cite {kowalski} or Remark 1 of \cite {nguyen1}.
\\
\\ Considering instead the average
\[ \sideset {}{'}\sum _{a=1}^q\left |\sum _{n\leq x\atop {n\equiv a(q)}}d(n)-\text { main term }\right |\]
Banks, Heath-Brown and Shparlinski \cite {banks} showed that averaging over $a$ is indeed enough to get a power-saving error term when $q$ goes all the way up to $x$. With their method Nguyen \cite {nguyen1} has the nice result that we get cancellation up to $x^{2/3}$ for $k=3$. But for $k\geq 4$ there is nothing past $\sqrt x$ even for this average case (or indeed any other average), and the arguments in these papers give just $\sqrt x$. As a brief summary:
\\
\begin{center}
\begin{tabular}{||c c c||} 
 \hline
$k$ & individual modulus & on average over $a=1,...,q$ \\ [0.5ex] 
 \hline\hline
$2$ & $2/3$ (Selberg/Hooley 1950's) & $1$ (Banks/Heath-Brown/Shparlinski 2005)\\ 
 \hline
$3$ & $1/2+1/46$ (Fouvry/Kowalski/Michel 2015) & 2/3 (Nguyen 2021) \\ 
 \hline
$4$ & $1/2$ (Linnik 1950's) & $1/2$ (Linnik 1950's) \\ 
 \hline
$\geq 5$ & $<1/2$  & $<1/2$ \\ 
 \hline
\end{tabular}
\end{center}
In this note, we show that for $k=4$ we have equidistribution on average for prime $q\leq x^{4/7}$, and so pass the square root barrier.
\\
\begin {theorem}\label {t1}
Let
\begin {eqnarray*}
E_{a/q}(s)&=&\sum _{n=1}^\infty \frac {d_4(n)}{n^s}e\left (\frac {na}{q}\right )\hspace {10mm}f_{a/q}(x)=Res_{s=1}\left \{ \frac {E_{a/q}(s)\left ((2x)^s-x^s\right )}{s}\right \} 
\\ \hspace {0mm}\Delta (a/q)&=&\sum _{x<n\leq 2x}d_4(n)e\left (\frac {na}{q}\right )-f_{a/q}(x).
\end {eqnarray*}
Then for prime $q\leq x^{4/7}$
\begin {eqnarray*}
\sum _{a=1}^q\left |\Delta (a/q)\right |^2\ll x^{3/2+\epsilon }q^{7/8}.
\end {eqnarray*}
\end {theorem}
Let $c_q(n)$ be Ramanujan's sum
\begin {equation}\label {seren}
c_q(n)=\sideset {}{'}\sum _{a=1}^qe\left (\frac {an}{q}\right )\hspace {10mm}\text {and let }\hspace {10mm}\theta _n(q)=\frac {\mu (q/(q,n))}{\phi (q/(q,n))}.
\end {equation}
It is straight-forward to establish that
\[ \sum _{n=1\atop {d|n}}^\infty \frac {d_4(n)\chi (n/d)}{n^s}\]
is holomorphic past $s=1$ for non-principal Dirichlet characters so, after breaking the sum according to the residue $b$ mod $q$, then according to $d=(b,q)$, and then using Dirichlet characters, we get
\begin {eqnarray*}
f_{a/q}(x)&=&
Res_{s=1}\left \{ \frac {(2x)^s-x^s}{s}\sum _{n=1}^\infty \frac {d_4(n)\theta _n(q/(q,a))}{n^s}\right \} =:F_x\left (\frac {q}{(q,a)}\right ).
\end {eqnarray*}
Then letting
\[ 
\mathcal M_x(q,a)=\frac {1}{q}\sum _{d|q}c_d(a)F_x(d)\hspace {10mm}E_x(q,a)=\sum _{x<n\leq 2x\atop {n\equiv a(q)}}d_4(n)-\mathcal M_x(q,a)\]
we get
\[ \sum _{N=1}^q\mathcal M_x(q,N)e\left (\frac {Nb}{q}\right )=F_x\left (\frac {q}{(q,b)}\right )\]
so
\begin {eqnarray*}
E_x(q,a)=\frac {1}{q}\sum _{b=1}^qe\left (-\frac {ab}{q}\right )\left (\sum _{x<n\leq 2x}d_4(n)e\left (\frac {nb}{q}\right )-\sum _{N=1}^q\mathcal M_x(q,N)e\left (\frac {Nb}{q}\right )\right )=\frac {1}{q}\sum _{b=1}^qe\left (-\frac {ab}{q}\right )\Delta (b/q)
\end {eqnarray*}
and therefore
\begin {eqnarray*}
\sum _{a=1}^q|E_{x}(q,a)|^2&=&\frac {1}{q}\sum _{b=1}^q\left |\Delta (b/q)\right |^2
\end {eqnarray*}
so Theorem \ref {t1} leads to
\begin {theorem}\label {t2}
For prime $q\leq x^{4/7}$
\begin {eqnarray*}
\sum _{a=1}^q|E_x(q,a)|\ll x^{3/4+\epsilon }q^{7/16}.
\end {eqnarray*}
\end {theorem}
\hspace {1mm}
\\ which in particular shows equidistribution on average up to $q\leq x^{4/7}$.
\\
\\ Our starting point will be Voronoi's summation formula for general $d_k(n)$, due to Ivi\` c:
\begin {lemma}\label {voronoi3}
Let $E_{h/q}(s)$ be as in Theorem \ref {t1}, let $w:[0,\infty )\rightarrow \mathbb R$ be smooth and of compact support, and let
\begin {eqnarray*}
\tilde \Delta (h/q)&=&\sum _{n=1}^\infty d_k(n)e\left (\frac {nh}{q}\right )w(n)-Res_{s=1}\left \{ E_{h/q}(s)\int _0^\infty w(t)t^{s-1}dt\right \} 
\\ U(X)&=&\frac {1}{2\pi i}\int _{(c)}\underbrace {\left (\frac {\Gamma (s/2)}{\Gamma ((1-s)/2)}\right )^k}_{\asymp 1/|t|^{k(1/2-\sigma )}}\frac {ds}{X^s}\hspace {10mm}\text {for }X>0\text { and }0<c<1/2-1/k
\\ N&=&\frac {\pi ^kn}{q^k}\hspace {20mm}\hat w_q(n)=\int _0^\infty w(t)U(Nt)dt
\\ R_{a_1,...,a_k}(h/q)&=&\sum _{x_1,...,x_k=1}^qe\left (\frac {a_1x_1+\cdot \cdot \cdot +a_kx_k-hx_1\cdot \cdot \cdot x_k}{q}\right )\hspace {10mm}A_{h/q}(n)=\frac {1}{2}\sum _{a_1\cdot \cdot \cdot a_k=n}R_{a_1,...,a_k}(h/q).
\end {eqnarray*}
Then for $(h,q)=1$
\begin {eqnarray*}
\tilde \Delta (h/q)&=&\frac {\pi ^{k/2}}{q^k}\sum _{n=1}^\infty A_{h/q}(n)\hat w_q(n)+\text {\emph {3 similar terms}}.
\end {eqnarray*}
\end {lemma}
\begin {proof}
This is Theorem 2 of \cite {ivic} (the definitions of the various quantities in (2.1), (2.6), (2.7), (3.2) and line -4 of page 213, and the conditions on $f$ at the bottom of page 212).
\end {proof}
With this lemma in place our argument for Theorem \ref {t1} is simple enough and we give a brief outline. We smooth and apply the above summation formula - the variance we want becomes something like
\[ \approx \frac {1}{q^{2k}}\sum _{n,m=1}^\infty \hat w_q(n)\overline {\hat w_q(m)}\sideset {}{'}\sum _{a=1}^qA_{a/q}(n)\overline {A_{a/q}(m)}+\mathcal O\left (Y^2\right ).\]
The $\hat w_q(n)$ transform can be made arbitrarily small once $n>T:=(xq/Y)^k/x$ and made $\ll x^{1/2-1/2k}/N^{1/2+1/2k}$ otherwise. Also, if we just focus on terms $(nm,q)=1$ then the $a$ sum may be calculated as $q^kc_q(n-m)$. So the above becomes
\[ \approx x^{1-1/k}q\sum _{n,m\ll T\atop {(nm,q)=1}}\frac {d_k(n)d_k(m)c_q(n-m)}{(nm)^{1/2+1/2k}}+Y^2\ll \left (\frac {xq}{Y}\right )^{k-1}q+Y^2\ll \left (x^{k-1}q^{k}\right )^{2/(k+1)}\]
which for $k=4$ gives us cancellation for $q$ up to $\sqrt x$ (and recovers the results of \cite {banks}, \cite {blomer}, \cite {nguyen1} for $k=2,3$). Consequently any saving would take us past the square root and we find this saving through the Petrow-Young bound in the evaluation of the sum of Ramanujan sums on the LHS. After Lemma \ref {ramanujan2} we briefly indicate how convexity alone only recovers the exponent $1/2$.
\\ 
\\ By adapting Lemmas \ref {ramanujan1} and \ref {ramanujan2} below results for some higher $k$ could be obtained as well (as well as for $k=3$). The restriction $q$ prime is hopefully not a serious issue, this being a common simplification in these kind of arguments. We return to these questions in future work, the point for now being simply a first step past 1/2.
\\
\\ We give some comments comparing the argument here with past work. The argument in \cite {banks} to bound $d(n)$ is based on average bounds for incomplete Kloosterman sums. Blomer \cite {blomer} then used an easier argument based on the Voronoi summation formula (also giving stronger results than \cite {banks}) and by separating a diagonal term Lau and Zhao \cite{lauzhao} improved this to an asymptotic formula. The argument outlined above is this argument, giving $2/3$ for $k=3$ as in \cite {nguyen1} and giving $1/2$ for $k=4$. Our improvement comes from a slight refinement - taking into account the oscillation due to the Ramanujan sums, which were bounded absolutely before, and where we can find the cancellation through the Petrow-Young bound.
\\
\\ For $k=2,3$ our refinement would immediately improve the $d_3(n)$ exponent in \cite {nguyen1} to $8/11$ as well as strengthen the errors in (\cite {banks}, \cite {blomer} and) \cite {lauzhao} (even without subconvexity), but the main point we want to make is going past $1/2$ for $d_4(n)$.
\\
\\ \begin {center}
\section {Lemmas}
\end {center}
Lemmas \ref {r4} and \ref {asum} deal with the $A_{h/q}(n)$ sums from the Ivi\' c/Voronoi formula. Then Lemmas \ref {w0} and \ref {w1} are concerned with bounding the transform $\hat w_q(n)$ in the formula. Finally we will need to use cancellation within sums of Ramanujan sums and this is Lemmas \ref {ramanujan1} and \ref {ramanujan2}. The proof of Theorem \ref {t1} then starts on page 9.
\\ \begin {lemma}\label {r4}
For $k=4$ let $A_{h/q}(n)$ be as in Lemma \ref {voronoi3} with $q$ prime and $(h,q)=1$. If $q|n$ then 
\[ A_{h/q}(n)\ll q^\epsilon (q^4,n)\]
and if $q|n$ but $q\nmid m$ then
\[ \sideset {}{'}\sum _{h=1}^qA_{h/q}(n)\overline {A_{h/q}(m)}\ll q^{1+\epsilon }(q^4,n).\]
\end {lemma}
\begin {proof}
We omit writing in $q^\epsilon $ factors and use $c_q(n)\ll (q,n)$ and $(a,q)=1\implies c_q(an)=c_q(n),c_q(a)=\mu (q)$; also for easier reading let's write down the definitions
\[ A_{h/q}(n)=\frac {1}{2}\sum _{abcd=n}\underbrace {\sum _{x,y,z,w=1}^qe\left (\frac {ax+by+cz+dw-hxyzw}{q}\right )}_{=R_{a,b,c,d}(h/q)}.\]
Suppose $q|n$ and that (say) $q|d$. Then
\begin {eqnarray}\label {goriad}
R_{a,b,c,d}(h/q)&=&\underbrace {q\sum _{x,y,z=1\atop {q|hxyz}}^qe\left (\frac {ax+by+cz}{q}\right )}_{=R_{a,b,c,d}(1/q)}=q\sum _{x,y=1\atop {(q,xy)|c}}^qe\left (\frac {ax+by}{q}\right )(q,xy).
\end {eqnarray}
If $(q,c)=1$ this sum is $c_q(a)c_q(b)$ and if $(q,c)>1$ it is
\begin {eqnarray*}
&=&\sum _{x,y=1}^{q-1}e\left (\frac {ax+by}{q}\right )+q\sum _{x=1}^{q-1}e\left (\frac {ax}{q}\right )+q\sum _{y=1}^{q-1}e\left (\frac {by}{q}\right )+q
=c_q(a)c_q(b)+q\left (c_q(a)+c_q(b)+1\right )
\end {eqnarray*}
so
\begin {eqnarray*}
R_{a,b,c,d}(h/q)\ll \underbrace {q}_{=(q,d)}\left (c_q(a)c_q(b)+\underbrace {\overbrace {q}^{=(q,c)}\left (c_q(a)+c_q(b)+1\right )}_{\text {if }(q,c)>1}\right )\notag 
\ll (q,a)(q,b)(q,c)(q,d)\ll (q^4,n)
\end {eqnarray*}
which gives the first case. For the second case \eqref {goriad} says
\begin {eqnarray*}
\sideset {}{'}\sum _{h=1}^qR_{a,b,c,d}(h/q)\overline {R_{a',b',c',d'}(h/q)}&=&\underbrace {R_{a,b,c,d}(1/q)}_{\ll (q^4,n)}\sum _{x,y,z,w=1}^qe\left (\frac {x+y+z+w}{q}\right )c_q(xyzw)
\end {eqnarray*}
and here the sum is
\begin {eqnarray*}
&=&\sideset {}{'}\sum _{x,y,z,w=1}^qe\left (\frac {x+y+z+w}{q}\right )+\phi (q)\sum _{x,y,z,w=1\atop {q|xyzw}}^qe\left (\frac {x+y+z+w}{q}\right )=\mu (q)^4+\phi (q)\mu (q)^3
\end {eqnarray*}
so the last $h$ sum is $\ll (q^4,n)q$ and we're done.
\end {proof}

\begin {lemma}\label {asum}
Suppose $q$ is prime and let $A_{h/q}(n)$ be as in Lemma \ref {voronoi3} with $k=4$. If $(nm,q)=1$ then
\begin {eqnarray*}
\sideset {}{'}\sum _{h=1}^qA_{h/q}(n)\overline {A_{h/q}(m)}&=&q^4d_4(n)d_4(m)c_q(n-m)+\mathcal O\left (q^{3+\epsilon }\right ).
\end {eqnarray*}
\end {lemma}
\begin {proof}
We have
\[ R_{a,b,c,d}(h/q)\overline {R_{a',b',c',d'}(h/q)}=\sum _{x,y,z,w=1\atop {x',y',z',w'=1}}^qe\left (\frac {x+y+z+w-x'-y'-z'-w'}{q}\right )c_q\left (xyzwm-x'y'z'w'n\right ).\]
Since
\begin {eqnarray*}
\sum _{w=1}^qc_q(wM-N)e\left (\frac {w}{q}\right )&=&qe\left (\frac {N\overline M}{q}\right )\left \{ \begin {array}{ll}1&\text { if }(q,M)=1\\ 0&\text { if not}\end {array}\right .
\end {eqnarray*}
we can perform the $w,w'$ summations to get that the above sum becomes
\[ q^2\sideset {}{'}\sum _{x,y,z=1\atop {x',y',z'=1\atop {xyzm\equiv x'y'z'n(q)}}}^qe\left (\frac {x+y+z-x'-y'-z'}{q}\right )=q^2\sideset {}{'}\sum _{x,y=1\atop {x',y'=1}}^qe\left (\frac {x+y-x'-y'}{q}\right )c_q\left (xym-x'y'n\right ).\]
Now using 
\begin {eqnarray*}
(NM,q)=1\implies \sideset {}{'}\sum _{y=1}^qc_q(yM-N)e\left (\frac {y}{q}\right )&=&qe\left (\frac {N\overline M}{q}\right )+1
\end {eqnarray*}
so that
\begin {eqnarray*}
\sideset {}{'}\sum _{y,y'=1}^qc_q(yM-y'N)e\left (\frac {y-y'}{q}\right )&=&qc_q\left (N-M\right )+\mu (q)
\end {eqnarray*}
the $\mathbf x,\mathbf y$ sum becomes
\[ \sideset {}{'}\sum _{x,x'=1}^qe\left (\frac {x-x'}{q}\right )\Big (qc_q\left (nx'-mx\right )+\mu (q)\Big )=q\Big (qc_q(n-m)+\mu (q)\Big )+\mu (q)^3\]
and the lemma follows.
\end {proof}
\begin {lemma}\label {w0}
Let $U(X)$ be as in Lemma \ref {voronoi3} with $k=4$. If $X\gg 1$ then
\begin {eqnarray*}
U(X)&=&\frac {e((8X)^{1/4})}{X^{3/8}}+\text {\emph { similar terms but of lower order}}
\end {eqnarray*}
and there is a polynomial $P$ such that for $X\leq 1$
\begin {eqnarray*}
U(X)&=&P(\log X).
\end {eqnarray*}
\end {lemma}
\begin {proof}
(This is Lemma 3 of \cite {ivic}.) We will need the integral representation
\begin {equation}\label {bessel}
\frac {1}{2\pi i}\int _{(c)}\frac {\Gamma (w)dw}{\Gamma (1+\nu -w)X^{2w}}=\frac {J_{\nu }(2X)}{X^\nu }\hspace {10mm}0<c<\frac {\nu }{2}+\frac {3}{4}.
\end {equation}
Let $c_n$ stand for absolute constants, different at each occurence. For $\sigma \geq 1/100$ and $t,A>0,B\geq 0$
\begin {eqnarray*}
\log \frac {\Gamma (As-B)}{\Gamma (1/2-(As-B))}&=&
\left (2(As-B)-1/2\right )\log s+\sum _{n\geq 0}\frac {c_n}{s^n}
\end {eqnarray*}
so
\begin {eqnarray*}
\log \underbrace {\frac {4^{8s-3}\Gamma (s)^4\Gamma (7/4-4s)}{\Gamma (1/2-s)^4\Gamma (4s-5/4)}}_{=:h(s)}&=&\log s+\sum _{n\geq 0}\frac {c_n}{s^n}
\hspace {5mm}\text { i.e. }\hspace {5mm}h(s)=\sum _{n\geq -1}\frac {c_n}{s^n}
\end {eqnarray*}
and therefore 
\begin {eqnarray}\label {capuccino}
U(X)&=&\int _{(1/9)}\frac {\Gamma (4s-5/4)h(s)ds}{4^{8s-3}\Gamma (7/4-4s)X^{2s}}\notag 
\\ &=&c_{-1}\int _{(1/9)}\frac {s\Gamma (4s-5/4)ds}{4^{8s-3}\Gamma (7/4-4s)X^{2s}}+\sum _{n\geq 0}c_n\int _{(1/9)}\frac {\Gamma (4s-5/4)ds}{s^n4^{8s-3}\Gamma (7/4-4s)X^{2s}}\notag 
\\ &=&\frac {c_{-1}}{X^{5/8}}\int _{(4/9-5/4)}\frac {w\Gamma (w)ds}{\Gamma (1/2-w)(4X^{1/4})^{2w}}+\frac {1}{X^{5/8}}\sum _{n\geq 0}c_n\int _{(4/9-5/4)}\frac {\Gamma (w)dw}{(w+5/4)^n\Gamma (1/2-w)(4X^{1/4})^{2w}}\notag 
\\ &=:&\frac {c_{-1}I(X)}{X^{5/8}}+\frac {1}{X^{5/8}}\sum _{n\geq 0}c_nI_n(X).
\end {eqnarray}
Then
\begin {eqnarray}\label {istanbul1}
I_0(X)&=&Res_{w=0}\left \{ \frac {\Gamma (w)dw}{\Gamma (1/2-w)(4X^{1/4})^{2w}}\right \} +\int _{(1/3)}\frac {\Gamma (w)dw}{\Gamma (1/2-w)(4X^{1/4})^{2w}}
\\ &=&\frac {c_0}{X^{1/2}}+\frac {J_{-1/2}(8X^{1/4})}{(4X^{1/4})^{-1/2}}\notag 
\end {eqnarray}
from \eqref {bessel}, and
\begin {eqnarray}\label {istanbul2}
I_1(X)&=&Res_{w=0}\left \{ \frac {\Gamma (w)dw}{(w+5/4)\Gamma (1/2-w)(4X^{1/4})^{2w}}\right \} +\int _{(1/3)}\frac {\Gamma (w)dw}{(w+5/4)\Gamma (1/2-w)(4X^{1/4})^{2w}}
\\ &=&\frac {c_n}{X^{1/2}}-\int _{(4c-5/4)}\left (1-\frac {7}{4(w+5/4)}\right )\frac {\Gamma (w)dw}{(1/2-w)\Gamma (1/2-w)(4X^{1/4})^{2w}}\notag 
\\ &=&\frac {c_n}{X^{1/2}}-\int _{(4c-5/4)}\frac {\Gamma (w)dw}{\Gamma (3/2-w)(4X^{1/4})^{2w}}-7\int _{(4c-5/4)}\frac {\Gamma (w)dw}{(w+5/4)\Gamma (3/2-w)(4X^{1/4})^{2w}}.\notag 
\end {eqnarray}
Except that 1/2 is replaced by 3/2, the first integral here is the same as that in \eqref {istanbul1} and the second is the same as that in \eqref {istanbul2}. So we can apply the same arguments to them to get
\begin {eqnarray*}
I_1(X)&=&\frac {c_1}{X^{1/2}}+\frac {J_{1/2}(8X^{1/4})}{(4X^{1/4})^{1/2}}+\text { similar}
\end {eqnarray*}
and the same arguments may be applied to all $I_n(X)$'s so that
\begin {eqnarray*}
\sum _{n\geq 0}c_nI_n(X)&=&\frac {c}{X^{1/2}}+\frac {J_{-1/2}(8X^{1/4})}{(4X^{1/4})^{-1/2}}+\text { similar.}
\end {eqnarray*}
Putting this with
\begin {eqnarray*}
\int _{(4/9-5/4)}\frac {w\Gamma (w)dw}{\Gamma (1/2-w)X^{w/2}}&=&\sum _{n=1}^\infty Res_{s=-n}\left \{ \frac {w\Gamma (w)}{\Gamma (1/2-w)X^{w/2}}\right \} =\frac {X^{1/4}}{\sqrt \pi }\sum _{n=1}^\infty \frac {(-1)^{n+1}(X^{1/4})^{2n-1}}{(2n-1)!}=\frac {X^{1/4}\sin (4X^{1/4})}{\sqrt \pi }
\end {eqnarray*}
in \eqref {capuccino} we get the first claim. For the second claim we just move the line of integration to the left to get 
\begin {eqnarray*}
U(X)&=&Res_{s=0}\left \{ \left (\frac {\Gamma (s/2)}{\Gamma \left ((1-s)/2\right )}\right )^4\frac {1}{X^s}\right \} .
\end {eqnarray*}
\end {proof}

\begin {lemma}\label {w1}
Take $1\leq Y\leq x$ and let $w:[0,\infty )\rightarrow \mathbb R$ be a smooth function satisfying
\begin {eqnarray*}
w(t)&=&\left \{ \begin {array}{ll}0&t\in [0,x-Y]\\ 1&t\in [x,2x]\\ 0&t\in [2x+Y,\infty )\end {array}\right .
\hspace {20mm}w^{(j)}(t)\ll \frac {1}{Y^j}\hspace {5mm}(j\geq 0)
\end {eqnarray*}
and let $N$, $U(X)$, $\hat w_q(n)$ be as in Lemma \ref {voronoi3} with $k=4$. 
If $N\gg 1/x$ then for any $j\in \mathbb N$
\begin {eqnarray*}
\hat w_q(n)&\ll &\frac {Y}{(Nx)^{3/8}}\left (\frac {x^3}{NY^4}\right )^{j/4}.
\end {eqnarray*}
\end {lemma}
\begin {proof}
Let $f(t)=8t^{1/4}$ and let
\[ M(X)=\frac {e(f(X))}{X^{3/8}}\hspace {5mm}\text { and }\hspace {5mm}\mathcal M(N)=\int _0^\infty w(t)M\left (Nt\right )dt.\]
Then
\begin {eqnarray*}
N^{3/8}\mathcal M(N)
&=&\int _0^\infty \frac {w(t)e(f(Nt))}{t^{3/8}}dt\hspace {10mm}(\star )
\\ &=&\frac {1}{4i\pi N^{1/4}}\int _0^\infty \frac {d}{dt}\left \{ w(t)t^{3/8}\right \} e(f(Nt))dt
\\ \text {so that }\hspace {10mm}N^{3/8}|\mathcal M(N)|&\leq &\frac {1}{4N^{1/4}}\left (\underbrace {\int _0^\infty \frac {w'(t)e(f(Nt))}{t^{3/8}}\cdot t^{3/4}dt}_{\text {like $(\star )$ but multiplied by }t^{3/4}}+\underbrace {\int _0^\infty \frac {w(t)e(f(Nt))}{t^{3/8}}\cdot \frac {dt}{t^{1/4}}}_{\text {like }(\star )\text { but multiplied by $1/t^{1/4}$}}\right ).
\end {eqnarray*}
If we do this $j$ times we get
\begin {eqnarray*}
N^{3/8}|\mathcal M(N)|&\leq &\frac {1}{4^jN^{j/4}}\underbrace {\sum }_{2^j\text { terms, each with $A+B=j$}}\left |\underbrace {\int _0^\infty \frac {w^{(A)}(t)e(f(Nt))}{t^{3/8}}\cdot \frac {t^{3A/4}}{t^{B/4}}dt}_{\ll x^{5/8-j/4}+Y^{1-A}x^{A-j/4-3/8}}\right |
\\ &\ll &\frac {Y}{x^{3/8}}\left (\frac {x^3}{NY^4}\right )^{j/4}
\end {eqnarray*}
which is the claim for $M(X)$ instead of $U(X)$, and now we're ok with Lemma \ref {w0}.
\end {proof}
\begin {lemma}\label {ramanujan1}
For $N,M\geq q$ with $q$ prime
\begin {eqnarray*}
\sum _{n\sim N\atop {m\sim M\atop {(nm,q)=1}}}d_4(n)d_4(m)c_q(n-m)&\ll &(NM)^\epsilon \left (\frac {NM}{q}+q^{5/4}\sqrt {NM}\left (\min \{ N,M\} ^{1/8}+\sqrt {q}\right ).
\right ).
\end {eqnarray*}
\end {lemma}
\begin {proof}
We will omit writing factors $(NM)^\epsilon $, suppose $N\leq M$, and take parameters $T$ to always be $\leq M^{10}$. We will make use of the bounds
\[ L_{\chi }(1/2+it)\ll (qt)^{1/6}\hspace {15mm}\int _{y}^{T}\frac {|L_\chi (1/2+iT)|^2ds}{t^{2/3}}\ll T^{1/3}+\frac {q^{1/3}}{y^{1/3}}+\frac {\sqrt q}{y^{2/3}}\hspace {15mm}\sideset {}{^*}\sum _{\chi }\int _{1}^{T}\frac {|L_{\chi }(s)|^4dt}{t}\ll q\]
which come from respectively Corollary 1.3 of \cite {petrowyoung}, Theorem\footnote {I can't access this, but the theorem is stated in e.g. Remark 3 of ``A note on the fourth moment of Dirichlet $L$-functions" from Bui and Heath-Brown} ... of \cite {motohashi} and Theorem 10.1 of \cite {topics}. The first two imply
\begin {eqnarray*}
\int _1^T\frac {|L_\chi (1/2+it)|^4dt}{t}&\ll &q^{2/3}\int _1^y\frac {dt}{t^{1/3}}+q^{1/3}\int _y^T\frac {|L_\chi (1/2+it)|^2dt}{t^{2/3}}
\ll q^{3/4}+(qT)^{1/3}
\end {eqnarray*}
and the first implies 
\[ \frac {1}{T}\int _{1/2}^{1}|L_\chi (\sigma +iT)|^4d\sigma \ll \frac {q^{2/3}}{T^{1/3}}\]
so for $\chi \not =\chi _0$ we get from Perron's formula
\begin {eqnarray*}
f(N):=\sum _{n\leq N}d_4(n)\chi (n)
&\ll &\frac {N}{T}+\sqrt N\left (\int _{1/2}^{1/2+iT}\frac {|L_\chi (s)|^4ds}{|s|}+\frac {q^{2/3}}{T^{1/3}}\right )\hspace {10mm}\text {(any $T\geq 1$)}
\\ &\ll &\sqrt N\left (\frac {\sqrt N}{T}+ q^{3/4}+(qT)^{1/3}\right )\notag 
\\ &\ll &N^{5/8}q^{1/4}+\sqrt Nq^{3/4}\notag 
\end {eqnarray*}
so, from the third quoted bound,
\begin {eqnarray}\label {ordiwadd}
\sum _{\chi \not =\chi _0}f(N)f(M)&\ll &\sum _{\chi \not =\chi _0}\left (\frac {M}{T}+\sqrt M\int _{1/2}^{1/2+iT}\frac {|L_\chi (s)|^4ds}{|s|}+\frac {q^{2/3}}{T^{1/3}}\right )\left (N^{5/8}q^{1/4}+\sqrt Nq^{3/4}\right )\notag 
\\ &\ll &q\sqrt M\left (N^{5/8}q^{1/4}+\sqrt Nq^{3/4}\right )\notag 
\\ &\ll &q^{5/4}\sqrt {NM}\left (N^{1/8}+\sqrt {q}\right ).
\end {eqnarray}
Since 
\begin {eqnarray*}
\sum _{n\leq N\atop {(n,q)=1}}d_4(n)c_q(n-m)&=&\sideset {}{'}\sum _{r=1}^qc_q(r-m)\sum _{n\leq N\atop {n\equiv r(q)}}d_4(n)
\\ &=&\frac {1}{\phi (q)}\sum _{\chi (q)}\left (\sideset {}{'}\sum _{r=1}^q\overline {\chi }(r)c_q(r-m)\right )f(N)
\\ &=&\frac {1}{\phi (q)}\left (\mu (q)c_q(m)+q\sum _{\chi \not =\chi _0(q)}\chi (m)\right )f(N)
\end {eqnarray*}
we get from \eqref {ordiwadd}
\begin {eqnarray*}
\sum _{n\leq N\atop {m\leq M\atop {(nm,q)=1}}}d_4(n)d_4(m)c_q(n-m)&\ll &\frac {|f(N)|}{\phi (q)}\sum _{m\leq M}d_4(m)|c_q(m)|+\frac {q}{\phi (q)}\sum _{\chi \not =\chi _0}|f(N)f(M)|
\\ &\ll &\frac {NM}{q}+q^{5/4}\sqrt {NM}\left (N^{1/8}+\sqrt {q}\right ).
\end {eqnarray*}
and we're done. 
\end {proof}

\begin {lemma}\label {ramanujan2}
Take $q,Q,T\leq x^{\mathcal O(1)}$ with $q$ prime, take $P$ any polynomial, and write $\Sigma _{n,m}^\sharp $ for a sum subject to $n,m>q$ and $(nm,q)=1$. Then up to an error $\ll x^\epsilon $
\begin {eqnarray*}
\sideset {}{^\sharp }\sum _{n,m\leq Q}d_4(n)d_4(m)c_q(n-m)P(\log n/Q)P(\log m/Q)&\ll &
\frac {Q^2}{q}+\sqrt qQ^{3/2}
\\ 
\sideset {}{^\sharp }\sum _{n\leq Q\atop {Q<m\leq T}}\frac {d_4(n)d_4(m)c_q(n-m)P(\log n/Q)}{m^{5/8}}e\left (\left (n/Q\right )^{1/4}\right )&\ll &\frac {(Tx)^{1/4}}{q}\left (\frac {QT^{3/8}}{q}+q^{5/4}\sqrt Q+q^{7/4}Q^{3/8}\right )
\\ 
\sideset {}{^\sharp }\sum _{Q<n,m\leq T}\frac {d_4(n)d_4(m)c_q(n-m)}{(nm)^{5/8}}e\left (\left (n/Q\right )^{1/4}-\left (m/Q\right )^{1/4}\right )&\ll &\frac {(Tx)^{1/2}}{q^2}\left (\frac {T^{3/4}}{q}+\frac {q^{5/4}}{Q^{1/8}}+\frac {q^{7/4}}{Q^{1/4}}\right ).
\end {eqnarray*}
\end {lemma}
\begin {proof}
We will omit writing $x^\epsilon $ factors and will write $K=1/4,d=5/8$. We have
\begin {align*}
&f(t)\ll 1\hspace {10mm}f'(t)\ll \frac {1}{t}&\text { for }f(t)=P(\log t/Q)&
\\ &f(t)\ll \frac {1}{t^d}\hspace {10mm}f'(t)\ll \frac {(t/Q)^K+1}{t^{d+1}}&\text { for }f(t)=\frac {e\left ((t/Q)^{1/4}\right )}{t^d}&
\end {align*}
so from Lemma \ref {ramanujan1}
\begin {eqnarray*}
&&
\sideset {}{^\sharp }\sum _{n\sim N\atop {m\sim M}}d_4(n)d_4(m)c_q(n-m)f(n)g(m)
\\ &&\hspace {10mm}\ll \hspace {4mm}\max _{t\sim N\atop {t'\sim M}}\Big (|f(t)g(t')|+N|f'(t)g(t)|+M|f(t)g'(t')|+NM|f'(t)g'(t')|\Big )
\\ &&\hspace {14mm}\times \hspace {4mm}\left (\frac {NM}{q}+q^{5/4}\sqrt {NM}\left (\min \{ N,M\} ^{1/8}+\sqrt {q}\right )\right )
\\ &&\hspace {10mm}\ll \hspace {4mm}\frac {NM}{q}+q^{5/4}\sqrt {NM}\left (\min \{ N,M\} ^{1/8}+\sqrt {q}\right )\hspace {35mm}\text {for }f,g=P(...)
\\ &&\hspace {10mm}\ll \hspace {4mm}\left ((M/Q)^{1/4}+1\right )\left (\frac {NM^{1-d}}{q}+q^{5/4}\sqrt NM^{1/2-d}\left (\min \{ N,M\} ^{1/8}+\sqrt {q}\right )\right )\hspace {10mm}\text {for }f,g=P(...),\frac {e(...)}{...}
\\ &&\hspace {10mm}\ll \hspace {4mm}\left ((N/Q)^{1/4}+1\right )\left ((M/Q)^{1/4}+1\right )
\\ &&\hspace {21mm}\times \hspace {4mm}\left (\frac {(NM)^{1-d}}{q}+q^{5/4}(NM)^{1/2-d}\min \{ N,M\} ^{1/8}+q^{7/4}(NM)^{1/2-d}\right )\hspace {5mm}\text {for }f,g=\frac {e(...)}{...}
\end {eqnarray*}
so the three sums of the claim are respectively
\begin {eqnarray*}
&\ll &\max _{n,m\leq Q}\left (\sideset {}{^\sharp }\sum _{n\sim N\atop {m\sim M}}d_4(n)d_4(m)c_q(n-m)f(n)g(m)\right )\ll \frac {Q^2}{q}+q^{5/4}Q^{9/8}+q^{7/4}Q
\\ &\ll &\max _{n\leq Q\atop {Q<m\leq T}}\left (\sideset {}{^\sharp }\sum _{n\sim N\atop {m\sim M}}d_4(n)d_4(m)c_q(n-m)f(n)g(m)\right )\ll \left (\frac {(Tx)^K}{q}+1\right )\left (\frac {QT^{1-d}}{q}+q^{5/4}Q^{9/8-d}+q^{7/4}Q^{1-d}\right )
\\ &\ll &\max _{Q<n,m\leq T}\left (\sideset {}{^\sharp }\sum _{n\sim N\atop {m\sim M}}d_4(n)d_4(m)c_q(n-m)f(n)g(m)\right )\ll \left (\frac {(Tx)^{K}}{q}+1\right )^2\left (\frac {T^{2-2d}}{q}+q^{5/4}Q^{1/2-d}+q^{7/4}Q^{1-2d}\right ).
\end {eqnarray*}

\end {proof}

We note that if we'd used instead the convexity bound $L_\chi (s)\ll (qt)^{1/4}$ in Lemma \ref {ramanujan1} we'd've found only
\[ \sum _{n\leq N}d_4(n)\chi (n)\ll q\sqrt N\]
leading to a total bound $\ll q^2\sqrt {NM}$ in that lemma. Using this input in Lemma \ref {ramanujan2} we'd get a bound containing at least the term $\ll (Tx)^{1/2}/Q^{1/4}$. In turn this would in \eqref {convexity} below give a bound
\[ \frac {x^{3}q^2}{Y^2}\]
leading to a total error $\ll x^{3/2}q$ for Theorem \ref {t1}, which is the usual exponent $1/2$ limit.

\begin {center}
\section {Proof of Theorem \ref {t1}}
\end {center}
Take $q\leq x$, let $E_{a/q}(s),f_{a/q}(x),\Delta (a/q)$ be as in Theorem \ref {t1}, and let $w(t),\hat w_q(n)$ be as in Lemmas \ref {w1} and \ref {voronoi3}. A
ssume that all bounds can include a factor $x^\epsilon $ which we don't write in explicitly. Let $q\leq Y\leq x$ be a parameter and write
\[ T:=x^{\epsilon -1}\left (\frac {xq}{Y}\right )^{4}
\]
so that, choosing a very large $j$, Lemma \ref {w1} says
\begin {eqnarray}\label {w2}
n>T\implies \hat w_q(n)&\ll &
\frac {x^{9}}{n^{5/4}}\left (\frac {x^3q^4}{nY^4}\right )^{(j-5)/4}
<\frac {x^{9-\epsilon (j-5)/4}}{n^{5/4}}\ll \frac {1}{x^{100}n^{5/4}}.
\end {eqnarray}
As in\footnote {specifically you need (2.2) and (2.3) to get expressions like (2.9)-(2.12), but involving four sums $\Sigma _{a,b,c,d}$ of expressions involving at most $(\log q)^3$ and $\gamma _0(a/q)\gamma _0(b/q)\gamma _0(c/q)$ terms, to then sum over $d$, and to bound the remaining sum $\Sigma _{q|abc}$ using $\gamma _0(a/q)\ll 1$; this is the only part of our paper where we don't give the details} \S 2 of \cite {ivic} the meromorphic part of $E_{a/q}(s)$ has coefficients $\ll 1/q$ so from
\begin {eqnarray*}
\text {`` the difference}\hspace {5mm}\frac {(2x)^s-x^s}{s}
-\int _0^\infty w(t)t^{s-1}dt\hspace {5mm}\text { is holomorphic and }\hspace {3mm}\ll Y|x^{s-1}|\hspace {2mm}\text {"}
\end {eqnarray*}
we get
\begin {eqnarray*}
\Delta (a/q)&=&\overbrace {\sum _{n=1}^\infty d_4(n)e\left (\frac {na}{q}\right )w(n)-Res_{s=1}\left \{ E_{a/q}(s)\int _0^\infty w(t)t^{s-1}dt\right \} }^{=:\tilde \Delta (a/q)}
\\ &&-\sum _{n\text { in two intervals}\atop {\text {of length $Y$}}}d_4(n)e\left (\frac {na}{q}\right )w(n)+\mathcal O\left (\frac {Y}{q}\right )
\end {eqnarray*}
which gives
\begin {eqnarray*}
\sideset {}{'}\sum _{a=1}^q|\Delta (a/q)|^2&\ll &
\sideset {}{'}\sum _{a=1}^q|\tilde \Delta (a/q)|^2+q\sum _{n,m\text { in two intervals}\atop {\text {of length $Y$}\atop {n\equiv m(q)}}}d_4(n)d_4(m)+\frac {Y^2}{q}\ll \sideset {}{'}\sum _{a=1}^q|\tilde \Delta (a/q)|^2+Y^2
\end {eqnarray*}
and so
\begin {eqnarray}\label {balwnsabloda}
\sum _{a=1}^q|\Delta (a/q)|^2=\sum _{d|q}\sideset {}{'}\sum _{a=1}^d|\Delta (a/d)|^2\ll \sum _{d|q}\left (\sideset {}{'}\sum _{a=1}^d|\tilde \Delta (a/d)|^2+Y^2\right ).
\end {eqnarray} 
Write $N=x/d^4$ and $Q=d^4/2x$. We will use Lemma \ref {w0} which says that for $n\leq Q$
\begin {eqnarray}\label {bound}
\text {for }n\leq Q\hspace {10mm}\hat w_d(n)&=&\int _0^\infty w(t)P\left (\log ((nt)^{1/4}/d)\right )dt\hspace {10mm}\ll x\notag 
\\ \text {for }Q\ll n\leq T\hspace {10mm}\hat w_d(n)&=&\int _0^\infty \frac {w(t)e\left ((8Nt)^{1/4}\right )dt}{(Nt)^{3/8}}+\underbrace {\text { lower order terms}}_{=:L}\notag 
\\ &=&\frac {1}{4\pi iN^{5/8}}\int _0^\infty \frac {d}{dt}\left \{ t^{3/8}w(t)\right \} e((8Nt)^{1/4})+L\notag 
\\ &=&\frac {d^{5/2}}{4\pi in^{5/8}}\int _0^\infty t^{3/8}w'(t)e\left (\frac {(8nt)^{1/4}}{q}\right )+L\notag 
\\ &\ll &\frac {x^{3/8}d^{5/2}}{n^{5/8}};
\end {eqnarray}
we will concentrate just on the main term of the last equality, the same argument clearly being applicable for the terms in $L$. Write $\Sigma _{n,m}^\sharp $ for a sum with $n,m>d$ and $(nm,d)=1$. From the first claim of Lemma \ref {ramanujan2}
\begin {eqnarray*}
&&\frac {1}{d^4}\sideset {}{^\sharp }\sum _{n,m\leq Q}d_4(n)d_4(m)c_d(n-m)\hat w_d(n)\overline {\hat w_d(m)}
\\ &&\hspace {10mm}=\hspace {4mm}\frac {1}{d^4}\int _0^\infty \int _0^\infty w(\mathbf t)\left (\sideset {}{^\sharp }\sum _{n,m\leq Q}d
_4(n)d_4(m)c_d(n-m)P\left (\log ((nt)^{1/4}/d)\right )P\left (\log ((mt')^{1/4}/d\right )\right )d\mathbf t
\\ &&\hspace {10mm}\ll \hspace {4mm}\frac {1}{d^4}\int _0^\infty \int _0^\infty |w(\mathbf t)|\left (\frac {Q^2}{d}+(dQ)^{5/4}+d^{7/4}Q\right )d\mathbf t
\\ &&\hspace {10mm}\ll \hspace {4mm}d^3+x^{3/4}d^{9/4}+xd^{7/4},
\end {eqnarray*}
from the second claim of Lemma \ref {ramanujan2}
\begin {eqnarray*}
&&\frac {1}{d^4}\sideset {}{^\sharp }\sum _{n\leq Q\atop {Q<m\leq T}}d_4(n)d_4(m)c_d(n-m)\hat w_d(n)\overline {\hat w_d(m)}
\\ &&\hspace {10mm}=\hspace {4mm}\frac {1}{d^{3/2}}\int _0^\infty \int _0^\infty t^{3/8}w(t)w'(t')\left (\sideset {}{^\sharp }\sum _{n\leq Q\atop {Q<m\leq T}}\frac {d_4(n)d_4(m)c_d(n-m)P\left (\log ((nt)^{1/4}/d)\right )}{m^{5/8}}e\left (-\frac {(8mt)^{1/4}}{d}\right )\right )d\mathbf t+L
\\ &&\hspace {10mm}\ll \hspace {4mm}\frac {1}{d^{3/2}}\int _0^\infty \int _0^\infty \frac {t^{3/8}|w(t)w'(t')|(Tx)^{1/4}}{d}\left (\frac {QT^{3/8}}{d}+d^{5/4}\sqrt Q+d^{7/4}Q^{3/8}\right )d\mathbf t
\\ &&\hspace {10mm}\ll \hspace {4mm}\left (\frac {x}{Y}\right )^{5/2}d^3+\frac {x^2d^{7/4}}{Y}
\end {eqnarray*}
and from the third claim of Lemma \ref {ramanujan2} 
\begin {eqnarray}\label {convexity}
&&\frac {1}{d^4}\sideset {}{^\sharp }\sum _{Q<n,m\leq T}d_4(n)d_4(m)c_d(n-m)\hat w_d(n)\overline {\hat w_d(m)}\notag 
\\ &&\hspace {10mm}=\hspace {4mm}d^{}\int _0^\infty \int _0^\infty \mathbf t^{3/8}w'(\mathbf t)\left (\sideset {}{^\sharp }\sum _{Q<n,m\leq T}\frac {d_4(n)d_4(m)c_d(n-m)}{(nm)^{5/8}}e\left (\frac {(8nt)^{1/4}-(8mt')^{1/4}}{d}\right )\right )d\mathbf t+L\hspace {10mm}
\\ &&\hspace {10mm}\ll \hspace {4mm}d^{}\int _0^\infty \int _0^\infty \frac {\mathbf t^{3/8}|w'(\mathbf t)|(Tx)^{1/2}}{d^2}\left (\frac {T^{3/4}}{d}+\frac {d^{5/4}}{Q^{1/8}}+\frac {d^{7/4}}{Q^{1/4}}\right )d\mathbf t\notag 
\\ &&\hspace {10mm}\ll \hspace {4mm}\left (\frac {x}{Y}\right )^5d^3+\left (\frac {x}{Y}\right )^2xd^{7/4};\notag 
\end {eqnarray}
also the bounds in \eqref {bound} give
\begin {eqnarray*}
\frac {1}{d^4}\sum _{n,m\ll T\atop {n\text { or }m\leq d\atop {(nm,d)=1}}}d_4(n)d_4(m)c_d(n-m)\hat w_d(n)\overline {\hat w_d(m)}&\ll &\frac {1}{d^4}\left (x^2\sum _{n\leq Q\atop {m\leq d}}+x^{11/8}d^{5/2}\sum _{Q<n\leq T\atop {m\leq d}}\frac {1}{n^{5/8}}\right )|c_d(n-m)|\ll \frac {x^{5/2}d}{Y^{3/2}}.
\end {eqnarray*}
The four bounds of the last sentence give together
\begin {equation}\label {lleuad}
\frac {1}{d^4}\sum _{n,m\leq T\atop {(nm,d)=1}}d_4(n)d_4(m)c_d(n-m)\hat w_d(n)\overline {\hat w_d(m)}\ll \left (\frac {x}{Y}\right )^5d^3+\left (\frac {x}{Y}\right )^2xd^{7/4}\hspace {2mm}=:\hspace {2mm}\mathcal E.
\end {equation}
From \eqref {bound}
\begin {eqnarray*}
\frac {1}{d^8}\sum _{n,m\leq T\atop {(nm,d)=1}}|\hat w_d(n)\overline {\hat w_d(m)}|\mathcal O(d^3)&\ll &\frac {1}{d^5}\left (x^2\sum _{n,m\leq Q}+x^{3/4}d^5\sum _{n,m\leq T}\frac {1}{(nm)^{5/8}}\right )\ll \left (\frac {x}{Y}\right )^3d^3
\end {eqnarray*}
so Lemma \ref {asum} and \eqref {lleuad} give
\begin {equation}\label {serenwib}
\frac {1}{d^{8}}\sum _{n,m\leq T\atop {(nm,d)=1}}\hat w_d(n)\overline {\hat w_d(m)}\sideset {}{'}\sum _{h=1}^dA_{h/d}(n)\overline {A_{h/d}(m)}\ll \mathcal E.
\end {equation}
From Lemma \ref {r4} 
\[ (nm,d)>1\implies \sideset {}{'}\sum _{h=1}^dA_{h/d}(n)\overline {A_{h/d}(m)}\ll d(d^4,n)(d^4,m)\]
so \eqref {bound} gives
\begin {eqnarray*}
&&\frac {1}{d^8}\sum _{n,m\leq T\atop {(nm,d)>1}}\hat w_d(n)\overline {\hat w_d(m)}\sideset {}{'}\sum _{h=1}^dA_{h/d}(n)\overline {A_{h/d}(m)}
\\ &&\hspace {10mm}\ll \hspace {4mm}\frac {1}{d^7}\left (x^2\sum _{n,m\ll Q}1+x^{3/4}d^{5}\sum _{n,m\ll T}\frac {1}{(nm)^{5/8}}\right )(d^4,n)(d^4,m)
\ll \left (\frac {x}{Y}\right )^3d^{}
\end {eqnarray*}
which with \eqref {serenwib} gives
\[ \frac {1}{d^{8}}\sum _{n,m\leq T}\hat w_d(n)\overline {\hat w_d(m)}\sideset {}{'}\sum _{h=1}^dA_{h/d}(n)\overline {A_{h/d}(m)}\ll \mathcal E.\]
From Lemma \ref {voronoi3}, \eqref {w2} and this we get
\begin {eqnarray*}
\frac {1}{\pi ^2}\sideset {}{'}\sum _{h=1}^d|\tilde \Delta (h/d)|^2&=&\frac {1}{d^{8}}\sum _{n,m\leq T}\hat w_d(n)\overline {\hat w_d(m)}\sideset {}{'}\sum _{h=1}^dA_{h/d}(n)\overline {A_{h/d}(m)}+\text { similar}+\mathcal O\left (\text { tiny }\right )
\hspace {2mm}\ll \hspace {2mm}\mathcal E
\end {eqnarray*}
so \eqref {balwnsabloda} gives
\begin {eqnarray*}
\sum _{a=1}^q|\Delta (a/q)|^2&\ll &\left (\frac {x}{Y}\right )^5q^3+\left (\frac {x}{Y}\right )^2xq^{7/4}+Y^2
\end {eqnarray*}
and Theorem \ref {t1} follows on optimising the second and third errors. 
\\

\begin {center}

\end {center}
\hspace {1mm}
\\
\\
\\
\\
\\ Contact details: 
\\ \emph {Tomos Parry
\\ Bilkent University, Ankara, Turkey
\\ tomos.parry1729@hotmail.co.uk}
\\
\\ On behalf of all authors, the corresponding author states that there is no conflict of interest. My manuscript has no associated data.

\end {document}